\documentclass{amsproc}%
\usepackage{graphicx}
\usepackage{amscd}
\usepackage{amsmath}
\usepackage{amsfonts}
\usepackage{amssymb}%
\theoremstyle{plain}

\numberwithin{equation}{section}

\begin{document}
\title[Determinants for $n\times n$ matrices]{Determinants for $n\times n$ matrices and the symmetric Newton formula in the
$3\times3$ case}
\author{J. Szigeti}
\address{Institute of Mathematics, University of Miskolc, Miskolc, Hungary 3515}
\email{matjeno@uni-miskolc.hu}
\author{L. van Wyk}
\address{Department of Mathematical Sciences, Stellenbosch University\\
P/Bag X1, Matieland 7602, Stellenbosch, South Africa }
\email{LvW@sun.ac.za}
\thanks{The first author was supported by OTKA K-101515 of Hungary and by the
TAMOP-4.2.1.B-10/2/KONV-2010-0001 project with support by the European Union,
co-financed by the European Social Fund.}
\thanks{The second author was supported by the National Research Foundation of South
Africa under Grant No.~UID 72375. Any opinion, findings and conclusions or
recommendations expressed in this material are those of the authors and
therefore the National Research Foundation does not accept any liability in
regard thereto.}
\thanks{The authors thank P. N. Anh, L. Marki and J. H. Meyer for fruitful consultations.}
\subjclass[2010]{ 15A15,15A24,15B33,16S50}
\keywords{$\mathbb{Z}_{2}$-graded, symmetric, left (right) determinants and
characteristic polynomials, Cayley-Hamilton identities, the symmetric
$3\times3$ Newton trace formula}

\begin{abstract}
One of the aims of this paper is to provide a short survey on the
$\mathbb{Z}_{2}$-graded, the symmetric and the left (right) generalizations of
the classical determinant theory for square matrices with entries in an
arbitrary (possibly non-commutative) ring. This will put us in a position to
give a motivation for our main results. We use the preadjoint matrix to
exhibit a general trace expression for the symmetric determinant. The
symmetric version of the classical Newton trace formula is also presented in
the $3\times3$ case.

\end{abstract}
\maketitle

\noindent1. INTRODUCTION

\bigskip

The universal notion of a determinant has a long history. Determinants of
matrices with entries in non-commutative rings have been considered by many
mathematicians, among them Cayley, Study, Ore, Dieudonn\'{e} and others. An
excellent recent survey dealing with almost all existing determinants is
[GGRW], where the so called Gelfand-Retakh quasideterminants are used as a
main organizing tool.

It seems that the symmetric determinant does not fit into the general
framework presented in [GGRW]. One of the aims of this paper is to provide a
short introduction to the symmetric and the corresponding left and right
versions of the classical determinant theory. The natural symmetrization of
the determinant formula and of the adjoint matrix lead to extremely useful
concepts. It turns out that these constructions can serve as a starting point
of a new symmetric determinant theory for square matrices over an arbitrary
ring. The most important feature of this theory is that it can be used to
solve systems of left (or right) linear equations and to exhibit left and
right Cayley-Hamilton identities for matrices over a Lie nilpotent ring (see
[Sz1, Sz2, Sz3, SzT]).

The algebra of $n\times n$ matrices over an exterior (Grassmann) algebra $E$
is our basic example. The Lie nilpotent property of $E$ plays a central role
in the various applications of the symmetric and the corresponding left and
right determinants.

The so called (and not widely known) $\mathbb{Z}_{2}$-graded determinant (in
[SSz]) is also designed to "attack" matrices over an exterior algebra. Its
construction heavily depends on the natural $\mathbb{Z}_{2}$-grading
$E=E_{0}\oplus E_{1}$. The $\mathbb{Z}_{2}$-graded determinant (and adjoint)
can be used to give an explicit inverse formula and to exhibit left and right
Cayley-Hamilton identities for an $n\times n$ matrix over $E$.

The "misterious" superdeterminant of a supermatrix due to Kantor and Trishin
(in [KT]) is another concept which is closely related to the $\mathbb{Z}_{2}%
$-grading of the exterior algebra. The treatment in [KT] leads to the solution
of certain special systems of left (or right) linear equations and to an
invariant Cayley-Hamilton identity for supermatrices. Unfortunately, the lack
of complete understanding prevents us from dealing with the KT-superdeterminant.

Using the fact that $E$ is a local ring, the Dieudonn\'{e} determinant is a
well defined $\cup_{n=1}^{\infty}\mathrm{GL}_{n}(E)\longrightarrow
E_{\mathrm{ab}}^{\times}$ map satisfying certain natural rules. Here
$E_{\mathrm{ab}}^{\times}$ denotes the Abelianized multiplicative group of
units in $E$\ (see [Ros]). The Dieudonn\'{e} determinant is an important tool
in algebraic K-theory, but we cannot use it to solve systems of linear
equations (over a local ring) and to derive Cayley-Hamilton identities.

In the rest of this introductory section we try to explain why we restrict our
attention to the matrix algebras $\mathrm{M}_{n}(E)$ and $\mathrm{M}_{n,t}(E)$.

The Cayley-Hamilton theorem and the corresponding trace identity play a
fundamental role in proving classical results about the polynomial and trace
identities of the $n\times n$ matrix algebra $\mathrm{M}_{n}(K)$ over a field
$K$ (see [Dr, DrF, Row1, Row2]). In case of $\mathrm{char}(K)=0$, Kemer's
pioneering work (see [Ke]) on the T-ideals of associative algebras revealed
the importance of the identities satisfied by the full $n\times n$ matrix
algebra $\mathrm{M}_{n}(E)$ and by the algebra of $(n,t)$ supermatrices
$\mathrm{M}_{n,t}(E)$, where%
\[
E=K\left\langle v_{1},v_{2},...,v_{i},...\mid v_{i}v_{j}+v_{j}v_{i}=0\text{
for all }1\leq i\leq j\right\rangle
\]
is the exterior (Grassmann) algebra generated by the infinite sequence of
anticommutative indeterminates $(v_{i})_{i\geq1}$.

Let $K\left\langle x_{1},x_{2},\ldots,x_{i},\ldots\right\rangle $ denote the
polynomial $K$-algebra generated by the infinite sequence $x_{1},x_{2}%
,\ldots,x_{i},\ldots$\ of non-commuting indeterminates. The prime T-ideals of
this (free associative $K$-)algebra\ are exactly the T-ideals of the
identities satisfied by $\mathrm{M}_{n}(K)$ for $n\geq1$. The T-prime T-ideals
are the prime T-ideals plus the T-ideals of the identities of $\mathrm{M}%
_{n}(E)$ for $n\geq1$ and of $\mathrm{M}_{n,t}(E)$ for $n-1\geq t\geq1$.
Another remarkable result is that for a sufficiently large $n\geq1$, any
T-ideal contains the T-ideal of the identities satisfied by $\mathrm{M}%
_{n}(E)$.

Accordingly, the importance of matrices (and supermatrices) over certain
non-commutative rings is an evidence in the theory of PI-rings, nevertheless
this fact has been obvious for a long time in other branches of algebra
(e.g.~in the structure theory of semisimple rings). Thus the algebras
$\mathrm{M}_{n}(E)$ and $\mathrm{M}_{n,t}(E)$ served as the main motivation
for the development of the symmetric and the $\mathbb{Z}_{2}$-graded determinants.

In Section 2 we deal with the Study determinant of a quaternionic matrix (see
[A, St]) and we show why a similar embedding approach does not work
for~$\mathrm{M}_{n}(E)$. The main result in Section 2 is based on an embedding
of the two generated exterior algebra $E^{(2)}$\ into a $2\times2$ matrix
algebra over a commutative ring and gives a Cayley-Hamilton identity of degree
$2n$\ in $\mathrm{M}_{n}(E^{(2)})$. Section 3 is devoted to a simplified
version of the $\mathbb{Z}_{2}$-graded determinant. In Section 4 we present a
short introduction to the symmetric theory of determinants, we collect and
explain some known results and point out some similarities and differences
between the (traditional) commutative base ring case and the general case.

We do not intend to give a detailed study of the relationships between the
symmetric and any of the already existing determinant notions. Nevertheless, a
thorough comparison between the symmetric and the $\mathbb{Z}_{2}$-graded
determinant would be essential. For a $2\times2$ matrix over $E$, the second
right (left) determinant is the double of the right (left) $\mathbb{Z}_{2}%
$-graded determinant. A similar comparison in the $3\times3$ case would
probably require computer calculations.

The treatment in Section 4 puts us in a position to give a motivation for our
new results\ in Sections 5 and 6. First we prove that $\mathrm{sdet}%
(A)=\mathrm{tr}(AA^{\ast})=\mathrm{tr}(A^{\ast}A)$, where $\mathrm{sdet}(A)$
is the symmetric determinant, $\mathrm{tr}(A)$ is the sum of the diagonal
entries and $A^{\ast}$ is the so called preadjoint matrix of the $n\times n$
matrix $A\in\mathrm{M}_{n}(R)$. Then we present the following symmetric
version of the Newton trace formula for a $3\times3$ matrix $A\in
\mathrm{M}_{3}(R)$:%
\[
\mathrm{sdet}(A)=\mathrm{tr}^{3}(A)-\mathrm{tr}(A)\cdot\mathrm{tr}%
(A^{2})-\mathrm{tr}(A\cdot\mathrm{tr}(A)\cdot A)-\mathrm{tr}(A^{2}%
)\cdot\mathrm{tr}(A)+\mathrm{tr}(A^{3})+\mathrm{tr}\left(  (A^{\top}%
)^{3}\right)  ,
\]
where $A^{\top}$ denotes the transpose of $A$. The symmetric characteristic
polynomial of this $A$\ and the corresponding general Cayley-Hamilton identity
are also presented by traces.

\bigskip

\noindent2. THE\ EMBEDDING\ OF\ $\mathbb{H}$\ AND\ THE\ STUDY\ DETERMINANT

\bigskip

There are well known embeddings of the complex number field $\mathbb{C}$ and
of the skew field $\mathbb{H}$\ of the real quaternions into matrices:%
\[
a+bi\longmapsto\left[
\begin{array}
[c]{cc}%
a & b\\
-b & a
\end{array}
\right]  , \ a+bi+cj+dk=a+bi+(c+di)j\longmapsto\left[
\begin{array}
[c]{cc}%
a+bi & c+di\\
-c+di & a-bi
\end{array}
\right]  .
\]
The above definitions provide injective $\mathbb{R}$-algebra homomorphisms
$\mu:\mathbb{C}\rightarrow\mathrm{M}_{2}(\mathbb{R})$ and $\nu:\mathbb{H}%
\rightarrow\mathrm{M}_{2}(\mathbb{C})$. There is a natural $\mu_{2}%
:\mathrm{M}_{2}(\mathbb{C)}\rightarrow\mathrm{M}_{2}(\mathrm{M}_{2}%
(\mathbb{R}))$ extension of $\mu$:%
\[
\mu_{2}\left(  \left[
\begin{array}
[c]{cc}%
z_{1,1} & z_{1,2}\\
z_{2,1} & z_{2,2}%
\end{array}
\right]  \right)  =\left[
\begin{array}
[c]{cc}%
\mu(z_{1,1}) & \mu(z_{1,2})\\
\mu(z_{2,1}) & \mu(z_{2,2})
\end{array}
\right]  .
\]
Since $\mathrm{M}_{2}(\mathrm{M}_{2}(\mathbb{R}))\cong\mathrm{M}%
_{4}(\mathbb{R})$, the composition $\vartheta=\mu_{2}\circ\nu:\mathbb{H}%
\longrightarrow\mathrm{M}_{4}(\mathbb{R})$ is the following map:%
\[
a+bi+cj+dk\longmapsto\left[
\begin{array}
[c]{cccc}%
a & b & c & d\\
-b & a & -d & c\\
-c & d & a & -b\\
-d & -c & b & a
\end{array}
\right]  .
\]
Using the natural extensions%
\[
\nu_{n}:\mathrm{M}_{n}(\mathbb{H)}\longrightarrow\mathrm{M}_{n}(\mathrm{M}%
_{2}(\mathbb{C}))\cong\mathrm{M}_{2n}(\mathbb{C})\text{ and }\vartheta
_{n}:\mathrm{M}_{n}(\mathbb{H)}\longrightarrow\mathrm{M}_{n}(\mathrm{M}%
_{4}(\mathbb{R}))\cong\mathrm{M}_{4n}(\mathbb{R}),
\]
an $n\times n$ matrix over $\mathbb{H}$ can be viewed as a $2n\times2n$ matrix
over $\mathbb{C}$ or as a $4n\times4n$ matrix over $\mathbb{R}$. Now we can
define the Study determinant of a quaternionic matrix $A\in\mathrm{M}%
_{n}(\mathbb{H})$ as the ordinary determinant $\mathrm{S}\det(A)=\det
_{\mathbb{C}}\nu_{n}(A)$ in~$\mathrm{M}_{2n}(\mathbb{C})$. If we take the
absolute value of the complex number $\det_{\mathbb{C}}\nu_{n}(A)$, then
$\left\vert \det_{\mathbb{C}}\nu_{n}(A)\right\vert ^{2}=\det_{\mathbb{R}%
}\vartheta_{n}(A)$. The Study determinant has some nice properties and it is
used frequently in differential geometry and Lie theory. The Cayley-Hamilton
identity for $\vartheta_{n}(A)$ yields the same identity (with real
coefficients) of degree $4n$ for~$A$ itself.

A similar approach to get a useful determinant notion and a Cayley-Hamilton
identity in $\mathrm{M}_{n}(E)$ is impossible. The reason is that the
infinitely generated exterior algebra $E$ cannot be embedded into a full
matrix algebra over a commutative ring ($E$ does not satisfy any of the
standard identities). On the other hand the embedding approach gives the following:

\bigskip

\noindent\textbf{Theorem 2.1.}\textit{ Let\ }$A\in\mathrm{M}_{n}(E^{(2)}%
)$\textit{ be an }$n\times n$\textit{ matrix over the two generated exterior
algebra}%
\[
E^{(2)}=K\left\langle v_{1},v_{2}\mid v_{1}v_{2}=-v_{2}v_{1},v_{1}^{2}%
=v_{2}^{2}=0\right\rangle .
\]
\textit{Then }$A$\textit{ satisfies a Cayley-Hamilton identity of the form}%
\[
A^{2n}+c_{2n-1}A^{2n-1}+\cdots+c_{1}A+c_{0}I=0
\]
\textit{where }$c_{i}\in K$\textit{, }$0\leq i\leq2n-1$\textit{.}

\bigskip

\noindent\textbf{Proof.} The assignments%
\[
1\longmapsto\left[
\begin{array}
[c]{cc}%
1 & 0\\
0 & 1
\end{array}
\right]  , \ v_{1}\longmapsto\left[
\begin{array}
[c]{cc}%
x & x\\
0 & -x
\end{array}
\right]  , \ v_{2}\longmapsto\left[
\begin{array}
[c]{cc}%
y & 0\\
-2y & -y
\end{array}
\right]
\]
define a $K$-embedding $\varepsilon:E^{(2)}\longrightarrow\mathrm{M}%
_{2}(K\left[  x,y\right]  /(x^{2},y^{2}))$, where $(x^{2},y^{2}%
)\vartriangleleft K\left[  x,y\right]  $ is the ideal generated by the
monomials $x^{2}$ and $y^{2}$. Now consider the induced $K$-embedding%
\[
\varepsilon_{n}:\mathrm{M}_{n}(E^{(2)})\longrightarrow\mathrm{M}%
_{n}(\mathrm{M}_{2}(K\left[  x,y\right]  /(x^{2},y^{2})))\cong\mathrm{M}%
_{2n}(K\left[  x,y\right]  /(x^{2},y^{2})).
\]
The trace of any $2\times2$ block
\[
\varepsilon(b_{0}+b_{1}v_{1}+b_{2}v_{2}+b_{3}v_{1}v_{2})=
\]%
\[
\left[
\begin{array}
[c]{cc}%
b_{0}+b_{1}x+b_{2}y-b_{3}xy+(x^{2},y^{2}) & b_{1}x-b_{3}xy+(x^{2},y^{2})\\
-2b_{2}y+2b_{3}xy+(x^{2},y^{2}) & b_{0}-b_{1}x-b_{2}y+b_{3}xy+(x^{2},y^{2})
\end{array}
\right]
\]
in $\varepsilon_{n}(A)$ is of the form $2b_{0}+(x^{2},y^{2})$. Since the trace
of $\varepsilon_{n}(A)$ is the sum of the traces of the diagonal $2\times2$
blocks, we have $\mathrm{tr}(\varepsilon_{n}(A))=2b+(x^{2},y^{2})$ for some
$b\in K$. The coefficients of the characteristic polynomial of $\varepsilon
_{n}(A)$ are rational polynomial expressions (with zero constant terms) of the
traces $\mathrm{tr}(\left(  \varepsilon_{n}(A)\right)  ^{k})=\mathrm{tr}%
(\varepsilon_{n}(A^{k}))$, $k\geq1$ (Newton formulae). Thus the
Cayley-Hamilton identity for $\varepsilon_{n}(A)$ is of the form%
\[
\left(  \varepsilon_{n}(A)\right)  ^{2n}+\left(  c_{2n-1}+(x^{2}%
,y^{2})\right)  \left(  \varepsilon_{n}(A)\right)  ^{2n-1}+\cdots
\]
\[
+\left(  c_{1}+(x^{2},y^{2})\right)  \varepsilon_{n}(A)+\left(  c_{0}%
+(x^{2},y^{2})\right)  I=0
\]
with $c_{i}\in K$, $0\leq i\leq2n-1$. It follows that%
\[
\varepsilon_{n}(A^{2n}+c_{2n-1}A^{2n-1}+\cdots+c_{1}A+c_{0}I)=
\]
\[
\varepsilon_{n}(A^{2n})+c_{2n-1}\varepsilon_{n}(A^{2n-1})+\cdots
+c_{1}\varepsilon_{n}(A)+c_{0}I=0
\]
holds in $\mathrm{M}_{2n}(K\left[  x,y\right]  /(x^{2},y^{2}))$, and thus the
injectivity of $\varepsilon_{n}$ gives the desired identity. $\square$

\bigskip

\noindent3. THE $\mathbb{Z}_{2}$-GRADED\ DETERMINANT

\bigskip

A $\mathbb{Z}_{2}$-grading of an (associative) ring $R$ is a pair
$(R_{0},R_{1})$, where $R_{0}$ and $R_{1}$ are additive subgroups of $R$ such
that $R=R_{0}\oplus R_{1}$ and $R_{i}R_{j}\subseteq R_{i+j}$ for all
$i,j\in\{0,1\}$ and $i+j$ is taken modulo $2$. The relation $R_{0}%
R_{0}\subseteq R_{0}$ ensures that $R_{0}$ is a subring of $R$. It is easy to
see that the existence of $1\in R$ implies that $1\in R_{0}$.

A $\mathbb{Z}_{2}$-grading $(R_{0},R_{1})$\ of the ring $R$ is called central
if $R_{0}\subseteq Z(R)$ (here $\mathrm{Z}(R)$ denotes the centre of $R$). The
condition $R_{0}\subseteq Z(R)$ implies the Lie nilpotence (of index $2$) of
$R$. The general notion of the $\mathbb{Z}_{2}$-graded determinant (in [SSz])
is defined for an $n\times n$ matrix over an arbitrary ring $R$ with a central
$\mathbb{Z}_{2}$-grading $(R_{0},R_{1})$.

In order to present a more natural and understandable treatment of the
$\mathbb{Z}_{2}$-graded determinant, we restrict ourselves to the case of the
well known central $\mathbb{Z}_{2}$-grading $E=E_{0}\oplus E_{1}$ of the
(infinitely generated) exterior algebra. If we add one more (anticommutative)
generator $w$ to the infinite sequence $(v_{i})_{i\geq1}$, we obtain an
extended%
\[
E_{w}=K\left\langle w,v_{1},v_{2},...,v_{i},...\mid w^{2}=wv_{j}+v_{j}%
w=v_{i}v_{j}+v_{j}v_{i}=0\text{\ for all }1\leq i\leq j\right\rangle
\]
exterior algebra. An $n\times n$ matrix $A\in\mathrm{M}_{n}(E)$ can be
uniquely written as $A=A_{0}+A_{1}$ with $A_{0}\in\mathrm{M}_{n}(E_{0})$ and
$A_{1}\in\mathrm{M}_{n}(E_{1})$. The companion matrix of $A$\ is defined as
$A_{0}+A_{1}w\in\mathrm{M}_{n}(E_{w}(0))$, where the notations $E_{w}%
(0)=(E_{w})_{0}$ and $E_{w}(1)=(E_{w})_{1}$ are used for the even and the odd
part of the (central) $\mathbb{Z}_{2}$-grading $E_{w}=(E_{w})_{0}\oplus
(E_{w})_{1}$ of $E_{w}$.

Since $E_{w}(0)$ is commutative, the ordinary determinant and the ordinary
adjoint of $A_{0}+A_{1}w$ are defined and can be written as
\[
\mathrm{\det}(A_{0}+A_{1}w)=d_{0}+d_{1}w\in E_{w}(0)\text{ and }%
\mathrm{adj}(A_{0}+A_{1}w)=B_{0}+B_{1}w\in\mathrm{M}_{n}(E_{w}(0)),
\]
where $d_{0}\in E_{0}$, $d_{1}\in E_{1}$, $B_{0}\in\mathrm{M}_{n}(E_{0})$,
$B_{1}\in\mathrm{M}_{n}(E_{1})$ and each of these objects is uniquely
determined by $A$. Clearly, $d_{0}=\det(A_{0})$, $B_{0}=\mathrm{adj}(A_{0})$
and the elements $d_{1},b_{i,j}^{(1)}\in E_{1}$ are also polynomial
expressions of the entries $a_{i,j}^{(0)}$ and $a_{i,j}^{(1)}$ (note that
$A_{0}=[a_{i,j}^{(0)}]$, $A_{1}=[a_{i,j}^{(1)}]$ and $B_{1}=[b_{i,j}^{(1)}]$).

\bigskip

\noindent\textbf{Theorem 3.1.}\textit{\ The elements of the product matrices}
\[
A(B_{0}+B_{1})=(A_{0}+A_{1})(B_{0}+B_{1})\mathit{\ }\text{\textit{and}%
}\mathit{\ }(B_{0}+B_{1})A=(B_{0}+B_{1})(A_{0}+A_{1})
\]
\textit{are contained in the subring }$E_{0}[d_{1}]$\textit{\ of }%
$E$\textit{\ generated by }$d_{1}$\textit{\ and the elements of }$E_{0}%
$\textit{, namely}:
\[
A(B_{0}+B_{1}),(B_{0}+B_{1})A\in\mathrm{M}_{n}(E_{0}[d_{1}]).
\]

\bigskip

The containment $E_{0}\subseteq\mathrm{Z}(E)$ implies that the subring
$E_{0}[d_{1}]\subseteq E$ is commutative (the elements of $E_{0}[d_{1}]$ are
polynomials of $d_{1}$ with coefficients in $E_{0}$). As a consequence of
Theorem 3.1 the determinant and the adjoint of the matrices $A(B_{0}%
+B_{1}),(B_{0}+B_{1})A\in\mathrm{M}_{n}(E_{0}[d_{1}])$ are defined. We call
\[
\mathrm{rgdet}(A)=\det(A(B_{0}+B_{1}))\text{ the right }\mathbb{Z}%
_{2}\text{-graded determinant}%
\]
and
\[
\mathrm{rgadj}(A)=(B_{0}+B_{1})\mathrm{adj}(A(B_{0}+B_{1}))\text{ the right
}\mathbb{Z}_{2}\text{-graded adjoint}%
\]
(with respect to $E=E_{0}\oplus E_{1}$) of the matrix $A\in\mathrm{M}_{n}(E)$.
Since
\[
A(B_{0}+B_{1})\mathrm{adj}(A(B_{0}+B_{1}))=\det(A(B_{0}+B_{1}))I
\]
in $\mathrm{M}_{n}(E_{0}[d_{1}])$, we immediately obtain (in $\mathrm{M}%
_{n}(E)$) that:
\[
A\mathrm{rgadj}(A)=\mathrm{rgdet}(A)I.
\]
\noindent\textbf{Proposition 3.2.} (i)\textit{\ If }$T\in\mathrm{GL}_{n}%
(E_{0})$\textit{\ is an invertible matrix and }$A\in\mathrm{M}_{n}%
(E)$\textit{, then }$\mathrm{rgdet}(TAT^{-1})=\mathrm{rgdet}(A)$\textit{\ and
}$\mathrm{rgadj}(TAT^{-1})=T(\mathrm{rgadj}(A))T^{-1}$\textit{.}

\noindent(ii)\textit{\ If }$A\in\mathrm{M}_{n}(E_{0})$\textit{, then
}$\mathrm{rgdet}(A)=(\det(A))^{n}$\textit{\ and }$\mathrm{rgadj}%
(A)=(\det(A))^{n-1}\mathrm{adj}(A)$\textit{.}

\bigskip

The polynomial ring $E[t]$\ inherits a natural (and central) $\mathbb{Z}_{2}%
$-grading $E[t]=E_{0}[t]\oplus E_{1}[t]$ from $E=E_{0}\oplus E_{1}$. We define
the right $\mathbb{Z}_{2}$-graded characteristic polynomial of a matrix
$A\in\mathrm{M}_{n}(E)$ as the right $\mathbb{Z}_{2}$-graded determinant (with
respect to $E[t]=E_{0}[t]\oplus E_{1}[t]$) of the matrix $tI-A\in
\mathrm{M}_{n}(E[t])$, where $I$ is the identity matrix in $\mathrm{M}_{n}%
(E)$:
\[
\chi_{A}(t)=\mathrm{rgdet}(tI-A)=\lambda_{0}+\lambda_{1}t+\cdots+\lambda
_{k}t^{k}\in E[t],\ \lambda_{0},\lambda_{1},...,\lambda_{k}\in E\text{ and
}\lambda_{k}\neq0.
\]
Since $\mathrm{GL}_{n}(E_{0})\subseteq\mathrm{GL}_{n}(E_{0}[t])$, an immediate
consequence of Proposition 3.2 is that $\chi_{TAT^{-1}}(t)=\chi_{A}(t)$ for
any invertible matrix $T\in\mathrm{GL}_{n}(E_{0})$.

\bigskip

\noindent\textbf{Proposition 3.3.}\textit{\ If }$\chi_{A}(t)=\lambda
_{0}+\lambda_{1}t+\cdots+\lambda_{k}t^{k}$\textit{\ is the right }%
$\mathbb{Z}_{2}$\textit{-graded characteristic polynomial of the }$n\times
n$\textit{\ matrix }$A\in\mathrm{M}_{n}(E)$\textit{, then }$k=n^{2}%
$\textit{\ and }$\lambda_{n^{2}}=1$\textit{, }$\lambda_{0}=\mathrm{rgdet}%
(-A)$\textit{.}

\bigskip

\noindent\textbf{Theorem 3.4.}\textit{\ If }$\chi_{A}(t)\in E[t]$\textit{\ is
the right }$\mathbb{Z}_{2}$\textit{-graded characteristic polynomial of an
}$n\times n$\textit{\ matrix }$A\in\mathrm{M}_{n}(E)$\textit{\ and }$h(t)\in
E[t]$\textit{\ is arbitrary, then the left substitution of }$A$\textit{\ into
the product polynomial }$\chi_{A}(t)h(t)=\mu_{0}+\mu_{1}t+\cdots+\mu_{m}t^{m}%
$\textit{\ is zero: }$I\mu_{0}+A\mu_{1}+\cdots+A^{m}\mu_{m}=0$\textit{.}

\bigskip

\noindent4. THE SYMMETRIC AND THE RIGHT (LEFT) DETERMINANTS

\bigskip

Let $\mathrm{S}_{n}$ denote the symmetric group of all permutations of the set
$\{1,2,\ldots,n\}$. For an $n\times n$ matrix $A=[a_{i,j}]$ over an arbitrary
(possibly non-commutative) ring or algebra $R$ with $1$, the element%
\[
\mathrm{sdet}(A)=\underset{\tau,\rho\in\mathrm{S}_{n}}{\sum}\mathrm{sgn}%
(\rho)a_{\tau(1),\rho(\tau(1))}\cdots a_{\tau(t),\rho(\tau(t))}\cdots
a_{\tau(n),\rho(\tau(n))}%
\]%
\[
=\underset{\alpha,\beta\in\mathrm{S}_{n}}{\sum}\mathrm{sgn}(\alpha
)\mathrm{sgn}(\beta)a_{\alpha(1),\beta(1)}\cdots a_{\alpha(t),\beta(t)}\cdots
a_{\alpha(n),\beta(n)}%
\]
of $R$ can be obviously considered as the symmetric determinant of $A$.

The preadjoint matrix $A^{\ast}=[a_{r,s}^{\ast}]$ of an $n\times n$ matrix
$A=[a_{i,j}]$ (over an arbitrary ring or algebra $R$ with $1$) is defined as
the following natural symmetrization of the classical adjoint:%
\[
a_{r,s}^{\ast}=\underset{\tau,\rho}{\sum}\mathrm{sgn}(\rho)a_{\tau
(1),\rho(\tau(1))}\cdots a_{\tau(s-1),\rho(\tau(s-1))}a_{\tau(s+1),\rho
(\tau(s+1))}\cdots a_{\tau(n),\rho(\tau(n))}%
\]%
\[
=\underset{\alpha,\beta}{\sum}\mathrm{sgn}(\alpha)\mathrm{sgn}(\beta
)a_{\alpha(1),\beta(1)}\cdots a_{\alpha(s-1),\beta(s-1)}a_{\alpha
(s+1),\beta(s+1)}\cdots a_{\alpha(n),\beta(n)}\text{ },
\]
where the first sum is taken over all $\tau,\rho\in\mathrm{S}_{n}$ with
$\tau(s)=s$ and $\rho(s)=r$ (the second sum is taken over all $\alpha,\beta
\in\mathrm{S}_{n}$ with $\alpha(s)=s$ and $\beta(s)=r$). We note that the
$(r,s)$ entry of $A^{\ast}$ is exactly the signed symmetric determinant
$(-1)^{r+s}\mathrm{sdet}(A_{s,r})$\ of the $(n-1)\times(n-1)$\ minor $A_{s,r}%
$\ of $A$ arising from the deletion of the $s$-th row and the $r$-th column of
$A$. If $R$\ is commutative, then $A^{\ast}=(n-1)!\mathrm{adj}(A)$, where
$\mathrm{adj}(A)$ denotes the ordinary adjoint of $A$.

The right adjoint sequence $(P_{k})_{k\geq1}$ of $A$ is defined by the
recursion: $P_{1}=A^{\ast}$ and $P_{k+1}=(AP_{1}\cdots P_{k})^{\ast}$ for
$k\geq1$. Originally the $k$-th right determinant was defined as the top left
entry of the product matrix $AP_{1}\cdots P_{k}$. These definitions were
introduced in [Sz1].

The above mentioned $k$-th right determinant is not invariant with respect to
the conjugate action of $\mathrm{GL}_{n}(\mathrm{Z}(R))$ on $\mathrm{M}%
_{n}(R)$. A more appropriate (and invariant) definition for the $k$-th right
determinant is the trace of $AP_{1}\cdots P_{k}$ (see [Do] and~[Sz3]):
\[
\mathrm{rdet}_{(k)}(A)=\mathrm{tr}(AP_{1}\cdots P_{k}).
\]
The left adjoint sequence $(Q_{k})_{k\geq1}$ can be defined analogously:
$Q_{1}=A^{\ast}$ and $Q_{k+1}=(Q_{k}\cdots Q_{1}A)^{\ast}$ for $k\geq1$. The
$k$-th left determinant of $A$ is%
\[
\mathrm{ldet}_{(k)}(A)=\mathrm{tr}(Q_{k}\cdots Q_{1}A).
\]
Note that $\mathrm{rdet}_{(k+1)}(A)=\mathrm{rdet}_{(k)}(AA^{\ast})$ and
$\mathrm{ldet}_{(k+1)}(A)=\mathrm{ldet}_{(k)}(A^{\ast}A)$. The basic
properties of these determinants are given in the following theorems.

\bigskip

\noindent\textbf{Theorem 4.1. }(see [Do], [Sz3])\textit{ If }$T\in
\mathrm{GL}_{n}(\mathrm{Z}(R))$\textit{ is an invertible matrix with entries
in the centre\ }$\mathrm{Z}(R)$\textit{ of }$R$\textit{, then}%
\[
\mathrm{tr}(T^{-1}AT)=\mathrm{tr}(A),\text{ }(T^{-1}AT)^{\ast}=T^{-1}A^{\ast
}T,
\]%
\[
\mathrm{rdet}_{(k)}(T^{-1}AT)=\mathrm{rdet}_{(k)}(A),\text{ }\mathrm{ldet}%
_{(k)}(T^{-1}AT)=\mathrm{ldet}_{(k)}(A).
\]

\bigskip

\noindent The next results shed light on the fact that we call $\mathrm{radj}%
_{(k)}(A)=nP_{1}\cdots P_{k}$ the $k$-th right adjoint and $\mathrm{ladj}%
_{(k)}(A)=nQ_{k}\cdots Q_{1}$ the $k$-th left adjoint of $A$.

\bigskip

\noindent\textbf{Theorem 4.2. }(see [Sz1], [Sz3])\textit{ The product matrices
}$A\mathrm{radj}_{(1)}(A)$\textit{ and }$\mathrm{ladj}_{(1)}(A)A$\textit{ in
}$\mathrm{M}_{n}(R)$\textit{\ can be written as}%
\[
A\mathrm{radj}_{(1)}(A)=nAA^{\ast}=\mathrm{tr}(AA^{\ast})I+C^{\prime
}=\mathrm{rdet}_{(1)}(A)I+C^{\prime}%
\]
\textit{and}%
\[
\mathrm{ladj}_{(1)}(A)A=nA^{\ast}A=\mathrm{tr}(A^{\ast}A)I+C^{\prime\prime
}=\mathrm{ldet}_{(1)}(A)I+C^{\prime\prime}\text{\textit{ ,}}%
\]
\textit{where }$I$\textit{ is the identity matrix,}\textrm{ }$\mathrm{tr}%
(C^{\prime})=\mathrm{tr}(C^{\prime\prime})=0$\textit{ and all entries of the
matrices }$C^{\prime}$\textit{ and }$C^{\prime\prime}$\textit{ are in the
additive commutator subgroup }$[R,R]$\textit{ of }$R$\textit{ generated by all
elements of the form }$[u,v]=uv-vu$\textit{, }$u,v\in R$\textit{.}

\bigskip

\noindent\textbf{Theorem 4.3. }(see [Sz1])\textit{ If the ring }$R$\textit{
satisfies the polynomial identity}%
\[
\lbrack\lbrack\lbrack\ldots\lbrack\lbrack x_{1},x_{2}],x_{3}],\ldots
],x_{k}],x_{k+1}]=0
\]
\textit{(}$R$\textit{ is Lie nilpotent of index }$k$\textit{), then the
products }$A\mathrm{radj}_{(k)}(A)$\textit{ and }$\mathrm{ladj}_{(k)}%
(A)A$\textit{ are scalar matrices in }$\mathrm{M}_{n}(R)$\textit{\ such that}%
\[
A\mathrm{radj}_{(k)}(A)=nAP_{1}\cdots P_{k}=\mathrm{rdet}_{(k)}%
(A)I,\mathrm{ladj}_{(k)}(A)A=nQ_{k}\cdots Q_{1}A=\mathrm{ldet}_{(k)}(A)I.
\]

\bigskip

\noindent If $R$\ is commutative, then $\mathrm{radj}_{(1)}(A)=\mathrm{ladj}%
_{(1)}(A)=nA^{\ast}=n!\mathrm{adj}(A)$ and%
\[
\mathrm{rdet}_{(k)}(A)=\mathrm{ldet}_{(k)}(A)=n\left\{  (n-1)!\right\}
^{1+n+n^{2}+\cdots+n^{k-1}}\left\{  \mathrm{\det}(A)\right\}  ^{n^{k-1}}.
\]
If $R$\ is Lie nilpotent of index $2$ and $\frac{1}{n}\in R$, then in Theorem 4.2

\noindent$C^{\prime}\in\mathrm{M}_{n}([R,R])\subseteq\mathrm{M}_{n}%
(\mathrm{Z}(R))$ implies that
\[
AA^{\ast}=\frac{1}{n}(\mathrm{tr}(AA^{\ast})I+C^{\prime})\in\mathrm{M}%
_{n}(\mathrm{Z}(R)[\mathrm{tr}(AA^{\ast})]),
\]
where the subring $\mathrm{Z}(R)[\mathrm{tr}(AA^{\ast})]$ generated by
$\mathrm{Z}(R)$ and $\mathrm{tr}(AA^{\ast})$ is commutative. Thus%
\[
\mathrm{rdet}_{(2)}(A)=\mathrm{rdet}_{(1)}(AA^{\ast})=n!\det(AA^{\ast}).
\]

Let $1\leq t\leq n-1$ be an integer and $R=R_{0}\oplus R_{1}$ be a
$\mathbb{Z}_{2}$-grading of $R$. Now $A\in\mathrm{M}_{n}(R)$ is called an
$(n,t)$ supermatrix if
\[
a_{i,j}\in R_{0}\text{ for all }1\leq i,j\leq t\text{ and }t+1\leq i,j\leq n,
\]
and%
\[
a_{i,j}\in R_{1}\text{ for all }1\leq i\leq t,t+1\leq j\leq n\text{ and
}t+1\leq i\leq n,1\leq j\leq t.
\]
Thus an $(n,t)$ supermatrix can be partitioned into square and rectangular
blocks as follows:%
\[
A=\left[
\begin{array}
[c]{cc}%
A_{1,1} & A_{1,2}\\
A_{2,1} & A_{2,2}%
\end{array}
\right]  ,
\]
where $A_{1,1}$ is a $t\times t$ and $A_{2,2}$ is an $(n-t)\times(n-t)$ square
matrix over $R_{0}$ and $A_{1,2}$ is a $t\times(n-t)$ and $A_{2,1}$ is an
$(n-t)\times t$ rectangular matrix over $R_{1}$. Clearly, the set of all
$(n,t)$ supermatrices $\mathrm{M}_{n,t}(R)$ is a subring (algebra) of
$\mathrm{M}_{n}(R)$.

\bigskip

\noindent\textbf{Theorem 4.4. }(see [Sz2])\textit{ If }$R=R_{0}\oplus R_{1}%
$\textit{ is a }$\mathbb{Z}_{2}$\textit{-grading of }$R$\textit{ and }%
$A\in\mathrm{M}_{n,t}(R)$\textit{, then }$A^{\ast}\in\mathrm{M}_{n,t}%
(R)$\textit{ and }$\mathrm{rdet}_{(k)}(A),\mathrm{ldet}_{(k)}(A)\in R_{0}%
$\textit{ for all }$1\leq k$\textit{.}

\bigskip

\noindent Let $R[z]$ denote the ring of polynomials of the single commuting
indeterminate $z$, with coefficients in $R$. The $k$-th right (left)
characteristic polynomial of $A$ is the $k$-th right (left) determinant of the
$n\times n$ matrix $zI-A$ in $\mathrm{M}_{n}(R[z])$:%
\[
p_{A,k}(z)=\mathrm{rdet}_{(k)}(zI-A)\text{ and }q_{A,k}(z)=\mathrm{ldet}%
_{(k)}(zI-A).
\]

\bigskip

\noindent\textbf{Theorem 4.5. }(see [Sz2])\textit{ If }$R=R_{0}\oplus R_{1}%
$\textit{ is a }$\mathbb{Z}_{2}$\textit{-grading of }$R$\textit{ and }%
$A\in\mathrm{M}_{n,t}(R)$\textit{, then }$p_{A,k}(z),q_{A,k}(z)\in R_{0}%
[z]$\textit{ for all }$1\leq k$\textit{.}

\bigskip

\noindent The above characteristic polynomials appear in the following
Cayley-Hamilton theorems.

\bigskip

\noindent\textbf{Theorem 4.6. }(see [Sz3])\textit{ The first right
characteristic polynomial }$p_{A,1}(z)\in R[z]$\textit{ of a matrix }%
$A\in\mathrm{M}_{n}(R)$\textit{\ is of the form}%
\[
p_{A,1}(z)=\lambda_{0}^{(1)}+\lambda_{1}^{(1)}z+\cdots+\lambda_{n-1}%
^{(1)}z^{n-1}+\lambda_{n}^{(1)}z^{n}%
\]
\textit{with }$\lambda_{0}^{(1)},\lambda_{1}^{(1)},\ldots,\lambda_{n-1}%
^{(1)},\lambda_{n}^{(1)}\in R$\textit{ and }$\lambda_{n}^{(1)}=n!$\textit{.
The product matrix}

\noindent$n(zI-A)(zI-A)^{\ast}$\textit{ can be written as}%
\[
n(zI-A)(zI-A)^{\ast}=p_{A,1}(z)I+C_{0}+C_{1}z+\cdots+C_{n}z^{n},
\]
\textit{where the matrices }$C_{i}\in\mathrm{M}_{n}(R)$\textit{ are uniquely
determined by }$A$\textit{, }$\mathrm{tr}(C_{i})=0$\textit{ and each entry of
}$C_{i}$\textit{ is in }$[R,R]$\textit{, i.e. }$C_{i}\in\mathrm{M}_{n}%
([R,R])$\textit{ for all }$0\leq i\leq n$. \textit{The right}%
\[
(\lambda_{0}^{(1)}I+C_{0})+A(\lambda_{1}^{(1)}I+C_{1})+\cdots+A^{n-1}%
(\lambda_{n-1}^{(1)}I+C_{n-1})+A^{n}(n!I+C_{n})=0
\]
\textit{and a similar left}%
\[
(\mu_{0}^{(1)}I+D_{0})+(\mu_{1}^{(1)}I+D_{1})A+\cdots+(\mu_{n-1}%
^{(1)}I+D_{n-1})A^{n-1}+(n!I+D_{n})A^{n}=0
\]
\textit{Cayley-Hamilton identity with right and left matrix coefficients hold
for }$A$\textit{.}

\bigskip

\noindent\textbf{Theorem 4.7. }(see [Sz1])\textit{ If the ring }$R$\textit{
satisfies the polynomial identity}%
\[
\lbrack\lbrack\lbrack\ldots\lbrack\lbrack x_{1},x_{2}],x_{3}],\ldots
],x_{k}],x_{k+1}]=0
\]
\textit{(}$R$\textit{ is Lie nilpotent of index }$k$\textit{), then the }%
$k$\textit{-th right characteristic polynomial}

\noindent$p_{A,k}(z)\in R[z]$\textit{ of a matrix }$A\in\mathrm{M}_{n}%
(R)$\textit{\ is of the form}%
\[
p_{A,k}(x)=\lambda_{0}^{(k)}+\lambda_{1}^{(k)}z+\cdots+\lambda_{n^{k}-1}%
^{(k)}z^{n^{k}-1}+\lambda_{n^{k}}^{(k)}z^{n^{k}},
\]
\textit{with }$\lambda_{0}^{(k)},\lambda_{1}^{(k)},\ldots,\lambda_{n^{k}%
-1}^{(k)},\lambda_{n^{k}}^{(k)}\in R$\textit{ and }$\lambda_{n^{k}}%
^{(k)}=n\left\{  (n-1)!\right\}  ^{1+n+n^{2}+\cdots+n^{k-1}}$\textit{.}

\noindent\textit{The right}%
\[
(A)p_{A,k}=I\lambda_{0}^{(k)}+A\lambda_{1}^{(k)}+\cdots+A^{n^{k}-1}%
\lambda_{n^{k}-1}^{(k)}+A^{n^{k}}\lambda_{n^{k}}^{(k)}=0
\]
\textit{and a similar left}%
\[
q_{A,k}(A)=\mu_{0}^{(k)}I+\mu_{1}^{(k)}A+\cdots+\mu_{n^{k}-1}^{(k)}A^{n^{k}%
-1}+\mu_{n^{k}}^{(k)}A^{n^{k}}=0
\]
\textit{Cayley-Hamilton identity with right and left scalar coefficients hold
for }$A$\textit{. We also have }$(A)u=v(A)=0$\textit{, where }$u(z)=p_{A,k}%
(z)h(z)$, $v(z)=h(z)q_{A,k}(z)$\textit{ and }$h(z)\in R[z]$\textit{ is
arbitrary.}

\bigskip

Now consider $E$ as a base ring and observe that $E$ is Lie nilpotent of index
$2$. Thus the above Theorems 4.3 and 4.7 apply to $\mathrm{M}_{n}(E)$. The
natural $\mathbb{Z}_{2}$-grading $E=E_{0}\oplus E_{1}$ allows us to apply
Theorems 4.4 and 4.5 to $\mathrm{M}_{n,t}(E)$. The most remarkable
consequences of these theorems are the following: $\mathrm{M}_{n}(E)$ is
integral over $E_{0}$ of degree $2n^{2}$ and $\mathrm{M}_{n,t}(E)$ is integral
over $E_{0}$ of degree $n^{2}$ (see [Sz1, Sz2]).

For a $2\times2$ matrix $A=[a_{i,j}]\in\mathrm{M}_{2}(E)$ we have
$A=A_{0}+A_{1}$ with

\noindent$A_{0}=[a_{i,j}^{(0)}]\in\mathrm{M}_{2}(E_{0})$ and $A_{1}%
=[a_{i,j}^{(1)}]\in\mathrm{M}_{2}(E_{1})$. Thus%
\[
\mathrm{adj}(A_{0}+A_{1}w)=\left[
\begin{array}
[c]{cc}%
a_{2,2}^{(0)}+a_{2,2}^{(1)}w & -a_{1,2}^{(0)}-a_{1,2}^{(1)}w\\
-a_{2,1}^{(0)}-a_{2,1}^{(1)}w & a_{1,1}^{(0)}+a_{1,1}^{(1)}w
\end{array}
\right]  =B_{0}+B_{1}w,
\]
whence%
\[
B_{0}+B_{1}=\left[
\begin{array}
[c]{cc}%
a_{2,2}^{(0)}+a_{2,2}^{(1)} & -a_{1,2}^{(0)}-a_{1,2}^{(1)}\\
-a_{2,1}^{(0)}-a_{2,1}^{(1)} & a_{1,1}^{(0)}+a_{1,1}^{(1)}%
\end{array}
\right]  =\left[
\begin{array}
[c]{cc}%
a_{2,2} & -a_{1,2}\\
-a_{2,1} & a_{1,1}%
\end{array}
\right]  =A^{\ast}%
\]
and%
\[
2\mathrm{rgdet}(A)=2\det(A(B_{0}+B_{1}))=2\det(AA^{\ast})=\mathrm{rdet}%
_{(1)}(AA^{\ast})=\mathrm{rdet}_{(2)}(A)
\]
follow. The comparison of $\mathrm{rgdet}(A)$ and $\mathrm{rdet}_{(2)}(A)$ for
a $3\times3$ matrix $A\in\mathrm{M}_{3}(E)$ is a challenging problem.

\bigskip

\noindent5. THE\ TRACE FORM\ OF\ THE\ SYMMETRIC\ DETERMINANT

\bigskip

If the base ring $R$ is commutative, then $\mathrm{tr}(AB)=\mathrm{tr}%
(BA)$\ for all $A,B\in\mathrm{M}_{n}(R)$. In spite of the fact that this well
known trace identity is no longer valid for matrices over a non-commutative
ring, the first left and first right determinants of $A$\ coincide (it was not
recognized in [Sz3]).

\bigskip

\noindent\textbf{Theorem 5.1.}\textit{ The traces of the product matrices
}$A^{\ast}A$\textit{ and }$AA^{\ast}$\textit{ are both equal to the symmetric
determinant of }$A$\textit{: }%
\[
\mathrm{rdet}_{(1)}(A)=\mathrm{tr}(AA^{\ast})=\mathrm{sdet}(A)=\mathrm{tr}%
(A^{\ast}A)=\mathrm{ldet}_{(1)}(A).
\]

\bigskip

\noindent\textbf{Proof.} We prove that $\mathrm{tr}(AA^{\ast})=\mathrm{sdet}%
(A)$. (The proof of $\mathrm{sdet}(A)=\mathrm{tr}(A^{\ast}A)$ is similar.) The
trace of a matrix is the sum of the diagonal entries, hence%
\[
\mathrm{tr}(A^{\ast}A)=\underset{1\leq r,s\leq n}{\sum}a_{r,s}^{\ast}a_{s,r}%
\]%
\[
=\underset{(\alpha,\beta,s)\in\Delta_{n}}{\sum}\mathrm{sgn}(\alpha
)\mathrm{sgn}(\beta)a_{\alpha(1),\beta(1)}\cdots a_{\alpha(s-1),\beta
(s-1)}a_{\alpha(s+1),\beta(s+1)}\cdots a_{\alpha(n),\beta(n)}a_{\alpha
(s),\beta(s)}%
\]%
\[
=\underset{\alpha^{\prime},\beta^{\prime}\in\mathrm{S}_{n}}{\sum}%
\mathrm{sgn}(\alpha^{\prime})\mathrm{sgn}(\beta^{\prime})a_{\alpha^{\prime
}(1),\beta^{\prime}(1)}\cdots a_{\alpha^{\prime}(t),\beta^{\prime}(t)}\cdots
a_{\alpha^{\prime}(n),\beta^{\prime}(n)}=\mathrm{sdet}(A),
\]
where $\Delta_{n}=\{(\alpha,\beta,s)\mid\alpha,\beta\in\mathrm{S}_{n},1\leq
s\leq n,\alpha(s)=s\}$ and the map

\noindent$(\alpha,\beta,s)\longmapsto(\alpha^{\prime},\beta^{\prime})$ with%
\[
\alpha^{\prime}=\left(
\begin{array}
[c]{ccccccc}%
1 & \ldots & s-1 & s & \ldots & n-1 & n\\
\alpha(1) & \ldots & \alpha(s-1) & \alpha(s+1) & \ldots & \alpha(n) & s
\end{array}
\right)
\]
and%
\[
\beta^{\prime}=\left(
\begin{array}
[c]{ccccccc}%
1 & \ldots & s-1 & s & \ldots & n-1 & n\\
\beta(1) & \ldots & \beta(s-1) & \beta(s+1) & \ldots & \beta(n) & \beta(s)
\end{array}
\right)
\]
is a $\Delta_{n}\longrightarrow\mathrm{S}_{n}\times\mathrm{S}_{n}$ bijection.
Since
\[
\mathrm{sgn}(\alpha^{\prime})=(-1)^{n-s}\mathrm{sgn}(\alpha)\text{ },\text{
}\mathrm{sgn}(\beta^{\prime})=(-1)^{n-s}\mathrm{sgn}(\beta)
\]
and%
\[
a_{\alpha(1),\beta(1)}\cdots a_{\alpha(s-1),\beta(s-1)}a_{\alpha
(s+1),\beta(s+1)}\cdots a_{\alpha(n),\beta(n)}a_{\alpha(s),\beta(s)}%
\]%
\[
=a_{\alpha^{\prime}(1),\beta^{\prime}(1)}\cdots a_{\alpha^{\prime}%
(t),\beta^{\prime}(t)}\cdots a_{\alpha^{\prime}(n),\beta^{\prime}(n)}\text{
},
\]
the proof is complete. $\square$

\bigskip

\noindent\textbf{Corollary 5.2.}\textit{ The first right and left
characteristic polynomials of a matrix }$A\in\mathrm{M}_{n}(R)$%
\textit{\ coincide: }$p_{A,1}(z)=q_{A,1}(z)$\textit{. Thus we have }%
$\lambda_{i}^{(1)}=\mu_{i}^{(1)}$\textit{ for all }$0\leq i\leq n$\textit{ in
the corresponding Cayley-Hamilton identities (see Theorem 4.6).}

\bigskip

In view of Theorem 5.1 and Corollary 5.2, for the above determinants and
characteristic polynomials, it is reasonable to use the terminology
"symmetric" instead of "first right" and "first left".

The following observation for $2\times2$\ matrices over the Grassmann algebra
is due to Domokos (see [Do]).

\bigskip

\noindent\textbf{Proposition 5.3.}\textit{ If }$A=[a_{i,j}]$\textit{ is in
}$\mathrm{M}_{2}(R)$\textit{, then}%
\[
\mathrm{rdet}_{(2)}(A)-\mathrm{ldet}_{(2)}(A)=\mathcal{S}_{4}(a_{1,1}%
,a_{1,2},a_{2,1},a_{2,2}),
\]
\textit{where }$\mathcal{S}_{4}(x_{1},x_{2},x_{3},x_{4})=\sum_{\sigma
\in\mathrm{S}_{4}}\mathrm{sgn}(\sigma)x_{\sigma(1)}x_{\sigma(2)}x_{\sigma
(3)}x_{\sigma(4)}$\textit{ is the standard polynomial of degree four.}

\bigskip

\noindent\textbf{Proof.} Using%
\[
A^{\ast}=\left[
\begin{array}
[c]{cc}%
a_{2,2} & -a_{1,2}\\
-a_{2,1} & a_{1,1}%
\end{array}
\right]
\]
and the products%
\[
AA^{\ast}=\left[
\begin{array}
[c]{cc}%
a_{1,1}a_{2,2}-a_{1,2}a_{2,1} & -a_{1,1}a_{1,2}+a_{1,2}a_{1,1}\\
a_{2,1}a_{2,2}-a_{2,2}a_{2,1} & -a_{2,1}a_{1,2}+a_{2,2}a_{1,1}%
\end{array}
\right]  ,
\]%
\[
A^{\ast}A=\left[
\begin{array}
[c]{cc}%
a_{2,2}a_{1,1}-a_{1,2}a_{2,1} & a_{2,2}a_{1,2}-a_{1,2}a_{2,2}\\
-a_{2,1}a_{1,1}+a_{1,1}a_{2,1} & -a_{2,1}a_{1,2}+a_{1,1}a_{2,2}%
\end{array}
\right]  ,
\]
a direct computation shows that%
\[
\mathrm{rdet}_{(2)}(A)-\mathrm{ldet}_{(2)}(A)=\mathrm{rdet}_{(1)}(AA^{\ast
})-\mathrm{ldet}_{(1)}(A^{\ast}A)=\mathrm{sdet}(AA^{\ast})-\mathrm{sdet}%
(A^{\ast}A)
\]%
\[
=(a_{1,1}a_{2,2}-a_{1,2}a_{2,1})(-a_{2,1}a_{1,2}+a_{2,2}a_{1,1})+(-a_{2,1}%
a_{1,2}+a_{2,2}a_{1,1})(a_{1,1}a_{2,2}-a_{1,2}a_{2,1})
\]%
\[
-(-a_{1,1}a_{1,2}+a_{1,2}a_{1,1})(a_{2,1}a_{2,2}-a_{2,2}a_{2,1})-(a_{2,1}%
a_{2,2}-a_{2,2}a_{2,1})(-a_{1,1}a_{1,2}+a_{1,2}a_{1,1})
\]%
\[
-(a_{2,2}a_{1,1}-a_{1,2}a_{2,1})(-a_{2,1}a_{1,2}+a_{1,1}a_{2,2})-(-a_{2,1}%
a_{1,2}+a_{1,1}a_{2,2})(a_{2,2}a_{1,1}-a_{1,2}a_{2,1})
\]%
\[
+(a_{2,2}a_{1,2}-a_{1,2}a_{2,2})(-a_{2,1}a_{1,1}+a_{1,1}a_{2,1})+(-a_{2,1}%
a_{1,1}+a_{1,1}a_{2,1})(a_{2,2}a_{1,2}-a_{1,2}a_{2,2})
\]%
\[
=\mathcal{S}_{4}(a_{1,1},a_{1,2},a_{2,1},a_{2,2}).\text{ }\square
\]

\bigskip

\noindent\textbf{Corollary 5.4.}\textit{ If }$A=[a_{i,j}]$\textit{ is in
}$\mathrm{M}_{2}(R)$\textit{, then}%
\[
p_{A,2}(z)-q_{A,2}(z)=\mathrm{rdet}_{(2)}(zI-A)-\mathrm{ldet}_{(2)}(zI-A)
\]%
\[
=\mathcal{S}_{4}(z-a_{1,1},-a_{1,2},-a_{2,1},z-a_{2,2})=\mathcal{S}%
_{4}(-a_{1,1},-a_{1,2},-a_{2,1},-a_{2,2})
\]%
\[
=\mathcal{S}_{4}(a_{1,1},a_{1,2},a_{2,1},a_{2,2})
\]
\textit{is a constant polynomial in }$R[z]$\textit{.}

\bigskip

\noindent6. THE\ SYMMETRIC NEWTON\ FORMULAE\ FOR\ $2\times2$ AND $3\times3$ MATRICES

\bigskip

If our base ring $R$ is commutative, then the well known Newton trace formulae
for $2\times2$ and $3\times3$ matrices are the following:%
\[
2\det(A)=\mathrm{tr}^{2}(A)-\mathrm{tr}(A^{2}),
\]%
\[
6\det(A)=\mathrm{tr}^{3}(A)-3\mathrm{tr}(A)\mathrm{tr}(A^{2})+2\mathrm{tr}%
(A^{3}).
\]

\bigskip

\noindent\textbf{Proposition 6.1.}\textit{ If }$R$\textit{ is an arbitrary
ring and }$A\in\mathrm{M}_{2}(R)$\textit{, then the symmetric analogue}%
\[
\mathrm{sdet}(A)=\mathrm{tr}^{2}(A)-\mathrm{tr}(A^{2})
\]
\textit{of the classical }$2\times2$\textit{ Newton formula holds. Notice that
}$\mathrm{sdet}(A)=2\det(A)$\textit{ in case of a commutative }$R$\textit{.}

\bigskip

\noindent\textbf{Proof.} Using%
\[
A=\left[
\begin{array}
[c]{cc}%
a & b\\
c & d
\end{array}
\right]  \text{ and }A^{2}=\left[
\begin{array}
[c]{cc}%
a^{2}+bc & ab+bd\\
ca+dc & cb+d^{2}%
\end{array}
\right]  ,
\]
we obtain that%
\[
\mathrm{tr}^{2}(A)-\mathrm{tr}(A^{2})=(a+d)^{2}-(a^{2}+bc+cb+d^{2}%
)=ad+da-bc-cb=\mathrm{sdet}(A).\text{ }\square
\]

\bigskip

\noindent\textbf{Theorem 6.2.}\textit{ If }$R$\textit{ is an arbitrary ring
and }$A\in\mathrm{M}_{3}(R)$\textit{, then the following symmetric analogue of
the classical }$3\times3$\textit{ Newton formula holds:}%
\[
\mathrm{sdet}(A)=\mathrm{tr}^{3}(A)-\mathrm{tr}(A)\cdot\mathrm{tr}%
(A^{2})-\mathrm{tr}(A\cdot\mathrm{tr}(A)\cdot A)-\mathrm{tr}(A^{2}%
)\cdot\mathrm{tr}(A)+\mathrm{tr}(A^{3})+\mathrm{tr}\left(  (A^{\top}%
)^{3}\right)  .
\]
\textit{Notice that}%
\[
\mathrm{sdet}(A)=6\det(A)\text{\textit{, }}\mathrm{tr}(A)\mathrm{tr}%
(A^{2})=\mathrm{tr}(A\cdot\mathrm{tr}(A)\cdot A)=\mathrm{tr}(A^{2}%
)\cdot\mathrm{tr}(A)\text{\textit{, }}\mathrm{tr}\left(  (A^{\top}%
)^{3}\right)  =\mathrm{tr}(A^{3})
\]
\textit{in case of a commutative }$R$\textit{.}

\bigskip

\noindent\textbf{Proof.} Using%
\[
A=\left[
\begin{array}
[c]{ccc}%
a & b & c\\
d & e & f\\
g & h & p
\end{array}
\right]  ,\text{ }A^{\top}=\left[
\begin{array}
[c]{ccc}%
a & d & g\\
b & e & h\\
c & f & p
\end{array}
\right]
\]%
\[
A^{2}=\left[
\begin{array}
[c]{ccc}%
a^{2}+bd+cg & ab+be+ch & ac+bf+cp\\
da+ed+fg & db+e^{2}+fh & dc+ef+fp\\
ga+hd+pg & gb+he+ph & gc+hf+p^{2}%
\end{array}
\right]
\]
and%
\[
A^{3}=\left[
\begin{array}
[c]{ccc}%
a^{2}+bd+cg & ab+be+ch & ac+bf+cp\\
da+ed+fg & db+e^{2}+fh & dc+ef+fp\\
ga+hd+pg & gb+he+ph & gc+hf+p^{2}%
\end{array}
\right]  \cdot\left[
\begin{array}
[c]{ccc}%
a & b & c\\
d & e & f\\
g & h & p
\end{array}
\right]  ,
\]
we obtain that%
\[
\mathrm{tr}(A^{2})=a^{2}+bd+cg+db+e^{2}+fh+gc+hf+p^{2}%
\]
and%
\[
\mathrm{tr}(A^{3})=(a^{2}+bd+cg)a+(ab+be+ch)d+(ac+bf+cp)g
\]%
\[
+(da+ed+fg)b+(db+e^{2}+fh)e+(dc+ef+fp)h
\]%
\[
+(ga+hd+pg)c+(gb+he+ph)f+(gc+hf+p^{2})p.
\]

We obtain a similar expression for $\text{\textrm{tr}}((A^{\top})^{3})$. Then
the proof can be completed by direct (but annoying) computation. $\square$

\bigskip

\noindent\textbf{Remark 6.3.} \textit{For a non-commutative ring }$R$\textit{,
the identity }$\mathrm{tr}\left(  (A^{\top})^{2}\right)  =\mathrm{tr}(A^{2}%
)$\textit{ holds for any }$A\in\mathrm{M}_{n}(R)$\textit{, but }%
$\mathrm{tr}\left(  (A^{\top})^{3}\right)  =\mathrm{tr}(A^{3})$\textit{ is not
valid even in the }$2\times2$\textit{\ case.}

\bigskip

\noindent\textbf{Theorem 6.4.}\textit{ If }$A$\textit{ is a }$3\times
3$\textit{ matrix over an arbitrary ring }$R$\textit{, then the symmetric
characteristic polynomial of }$A$\textit{ in }$R[z]$\textit{\ is}%
\[
p_{A,1}(z)=q_{A,1}(z)=\mathrm{sdet}(zI-A)=6z^{3}-6\mathrm{tr}(A)z^{2}%
+3(\mathrm{tr}^{2}(A)-\mathrm{tr}(A^{2}))z-\mathrm{sdet}(A).
\]

\bigskip

\noindent\textbf{Proof.} Using%
\[
\mathrm{tr}(zI-A)=3z-\mathrm{tr}(A),\mathrm{tr}^{3}(zI-A)=27z^{3}%
-27\mathrm{tr}(A)z^{2}+9\mathrm{tr}^{2}(A)z-\mathrm{tr}^{3}(A),
\]%
\[
(zI-A)^{2}=z^{2}I-zA-Az+A^{2},\text{ }\mathrm{tr}((zI-A)^{2})=3z^{2}%
-2\mathrm{tr}(A)z+\mathrm{tr}(A^{2}),
\]%
\[
\mathrm{tr}(zI-A)\cdot\mathrm{tr}((zI-A)^{2})=9z^{3}-9\mathrm{tr}%
(A)z^{2}+2\mathrm{tr}^{2}(A)z+3\mathrm{tr}(A^{2})z-\mathrm{tr}(A)\mathrm{tr}%
(A^{2}),
\]%
\[
\mathrm{tr}((zI-A)^{2})\cdot\mathrm{tr}(zI-A)=9z^{3}-9\mathrm{tr}%
(A)z^{2}+2\mathrm{tr}^{2}(A)z+3\mathrm{tr}(A^{2})z-\mathrm{tr}(A^{2}%
)\mathrm{tr}(A),
\]%
\[
\mathrm{tr}\left(  (zI-A)\cdot\mathrm{tr}(zI-A)\cdot(zI-A)\right)  =
\]%
\[
=\mathrm{tr}\left(  3Iz^{3}-3Az^{2}-\mathrm{tr}(A)Iz^{2}-3Az^{2}%
+\mathrm{tr}(A)Az+3A^{2}z+A\mathrm{tr}(A)z-A\mathrm{tr}(A)A\right)  =
\]%
\[
=9z^{3}-9\mathrm{tr}(A)z^{2}+2\mathrm{tr}^{2}(A)z+3\mathrm{tr}(A^{2}%
)z-\mathrm{tr}(A\cdot\mathrm{tr}(A)\cdot A),
\]%
\[
(zI-A)^{3}=z^{3}I-z^{2}A-zAz-Az^{2}+A^{2}z+AzA+zA^{2}-A^{3},
\]%
\[
\mathrm{tr}((zI-A)^{3})=3z^{3}-3\mathrm{tr}(A)z^{2}+3\mathrm{tr}%
(A^{2})z-\mathrm{tr}(A^{3}),
\]%
\[
\mathrm{tr}((zI-A)^{\top})^{3})=3z^{3}-3\mathrm{tr}(A^{\top})z^{2}%
+3\mathrm{tr}((A^{\top})^{2})z-\mathrm{tr}((A^{\top})^{3})=
\]%
\[
=3z^{3}-3\mathrm{tr}(A)z^{2}+3\mathrm{tr}(A^{2})z-\mathrm{tr}((A^{\top})^{3})
\]
and Theorem 6.2, we obtain that%

\[
\mathrm{sdet}(zI-A)=\mathrm{tr}^{3}(zI-A)
\]%
\[
-\!\mathrm{tr}(zI-A)\cdot\mathrm{tr}((zI-A)^{2})\!-\!\mathrm{tr}\left(
(zI\!-\!A)\cdot\mathrm{tr}(zI\!-\!A)\cdot(zI\!-\!A)\right)  \!-\!\mathrm{tr}%
((zI\!-\!A)^{2})\cdot\mathrm{tr}(zI\!-\!A)
\]%
\[
+\mathrm{tr}((zI-A)^{3})+\mathrm{tr}\left(  ((zI-A)^{\top})^{3}\right)
=6z^{3}-6\mathrm{tr}(A)z^{2}+3(\mathrm{tr}^{2}(A)-\mathrm{tr}(A^{2}%
))z-\mathrm{sdet}(A).\text{ }\square
\]

\bigskip

\noindent\textbf{Corollary 6.5. }\textit{If }$A\in\mathrm{M}_{3}(R)$\textit{,
then Theorem 4.6 gives the existence of }$3\times3$ \textit{matrices }%
$C_{i},D_{i}$\textit{ }$(0\leq i\leq3)$\textit{ with entries in }%
$[R,R]$\textit{ such that}%
\[
(-\mathrm{sdet}(A)I+C_{0})+(3(\mathrm{tr}^{2}(A)-\mathrm{tr}(A^{2}%
))I+C_{1})A+(-6\mathrm{tr}(A)I+C_{2})A^{2}+(6I+C_{3})A^{3}=0
\]
\textit{and}%
\[
(-\mathrm{sdet}(A)I+D_{0})+A(3(\mathrm{tr}^{2}(A)-\mathrm{tr}(A^{2}%
))I+D_{1})+A^{2}(-6\mathrm{tr}(A)I+D_{2})+A^{3}(6I+D_{3})=0.
\]

\bigskip

\noindent\textbf{Corollary 6.6. }\textit{If }$\frac{1}{6}\in R$\textit{ and
}$A\in\mathrm{M}_{3}(R)$\textit{ such that}%
\[
\mathrm{tr}(A)=\mathrm{tr}(A^{2})=\mathrm{tr}(A^{3})=\mathrm{tr}\left(
(A^{\top})^{3}\right)  =0,
\]
\textit{then }$\mathrm{sdet}(A)=0$\textit{ and}%
\[
A^{3}=C_{0}+C_{1}A+C_{2}A^{2}+C_{3}A^{3}=D_{0}+AD_{1}+A^{2}D_{2}+A^{3}D_{3}%
\]
\textit{for some }$3\times3$ \textit{matrices }$C_{i},D_{i}$\textit{ }$(0\leq
i\leq3)$\textit{ with entries in }$[R,R]$\textit{. Thus }$A^{3}\in
\mathrm{M}_{3}(T)$\textit{, where }$T=R[R,R]\cap\lbrack R,R]R$\textit{ is the
intersection of the left and right ideals }$R[R,R]$\textit{ and }%
$[R,R]R$\textit{ of }$R$\textit{.}

\bigskip

We close the paper by the following:

\bigskip

\noindent\textbf{Problem 6.7.}\textit{ If }$R$\textit{\ is a commutative ring,
then the Newton formula for a }$4\times4$\textit{\ matrix }$A\in\mathrm{M}%
_{4}(R)$\textit{ is}%
\[
24\det(A)=\mathrm{tr}^{4}(A)-6\mathrm{tr}^{2}(A)\mathrm{tr}(A^{2}%
)+3\mathrm{tr}^{2}(A^{2})+8\mathrm{tr}(A)\mathrm{tr}(A^{3})-6\mathrm{tr}%
(A^{4}).
\]
\textit{Find the symmetric analogue of the above formula for }$\mathrm{sdet}%
(A)$\textit{ over an arbitrary ring }$R$\textit{.}

\bigskip

\noindent\textbf{Acknowledgment.} The authors thank the referee for his/her
kind help to improve the exposition of the paper.

\bigskip

\noindent REFERENCES

\bigskip

\begin{enumerate}
\item {}[A] H. Aslaksen, \textit{Quaternionic determinants,}
Math.~Intelligencer 18 (3) (1996), 57--65.

\item {}[Do] M. Domokos, \textit{Cayley-Hamilton theorem for }$2\times
2$\textit{\ matrices over the Grassmann algebra}, J. Pure Appl. Algebra 133
(1998), 69-81.

\item {}[Dr] V. Drensky, \textit{Free Algebras and PI-Algebras},
Springer-Verlag, 2000.

\item {}[DrF] V. Drensky and E. Formanek, \textit{Polynomial Identity Rings},
Birkh\"{a}user-Verlag, 2004.

\item {}[GGRW] I. Gelfand, S. Gelfand, V. Retakh and R. L. Wilson,
\textit{Quasideterminants, }Adv.~Math. 193 (2005), 56--141.

\item {}[KT] I. Kantor and I. Trishin, \textit{On a concept of determinant in
the supercase,} Comm. Algebra~22 (10) (1994), 3679-3739.

\item {}[Ke] A. R. Kemer, \textit{Ideals of Identities of Associative
Algebras}, Translations of Math. Monographs, Vol. 87 (1991), AMS Providence,
Rhode Island.

\item {}[Ros] J. Rosenberg, \textit{Algebraic K-Theory and its Applications,}
GTM Springer, 1994.

\item {}[Row1] L. H. Rowen,\textit{\ Polynomial Identities in Ring Theory,}
Academic Press, New York, 1980.

\item {}[Row2] L. H. Rowen,\textit{\ Ring Theory Vol.~I, II,} Academic Press,
New York, 1988.

\item {}[SSz] S. Sehgal and J. Szigeti: \textit{Matrices over centrally
$\mathbb{Z}_{2}$-graded rings}, Beitr\"age Algebra Geom.~(Berlin) 43 (2)
(2002), 399-406.

\item {}[St] E. Study,\textit{ Zur Theorie der linearen Gleichungen,} Acta
Math.~42 (1920), 1-61.

\item {}[Sz1] J. Szigeti, \textit{New determinants and the Cayley-Hamilton
theorem for matrices over Lie nilpotent rings}, Proc.~Amer.~Math.~Soc.~125
(1997), 2245-2254.

\item {}[Sz2] J. Szigeti, \textit{On the characteristic polynomial of
supermatrices}, Israel J.~Math.~107 (1998), 229-235.

\item {}[Sz3] J. Szigeti, \textit{Cayley-Hamilton theorem for matrices over an
arbitrary ring}, Serdica Math.~J. 32 (2006), 269-276.

\item {}[SzT] J. Szigeti and Zs.~Tuza: \textit{Solving systems of linear
equations over Lie-nilpotent rings}, Linear and Multilinear Algebra~42 (1997), 43-51.
\end{enumerate}

\end{document}